\documentclass[reqno]{amsart}
\usepackage[cp1251]{inputenc}
\usepackage[english,russian]{babel}
\usepackage{amsmath}
\usepackage{amssymb}
\usepackage{amsfonts}
\usepackage{amsthm}

\def\negskp{\unskip\kern -2pt\unskip}
\newtheorem{conjecture}{Гипотеза}
\newtheorem{theorem}{Теорема}

\numberwithin{conjecture}{section}
\numberwithin{theorem}{section}
\numberwithin{lemma}{section}
\numberwithin{corollary}{section}
\numberwithin{proposition}{section}
\numberwithin{equation}{section}

\begin{document}
\thispagestyle{empty}
\title[О контрпримере Шарипова \dots]
{О контрпримере Шарипова к гипотезе Хабибуллина} 
\author{Р. А. Баладай}
\address{Рустам Алексеевич Баладай, 
\newline\hphantom{iii} ул. Дорожная 7, село Бураево, 
\newline\hphantom{iii} Республика Башкортостан, 
\newline\hphantom{iii} 452960, Россия}
\email{baladaichik@mail.ru}

\thanks{\sc Baladai R. A.
On Sharipov's counterexample to Habibullin's conjecture.}
\thanks{\copyright \ 2010 Баладай Р. А}
\maketitle 
{
\small
\begin{quote}
\noindent{\bf Аннотация. } Гипотеза Хабибуллина имеет три эквивалентные
формулировки. Недавно Шариповым был построен контрпример к этой гипотезе 
в одной из трех формулировок. В данной работе этот контрпример 
переносится на случай двух других формулировок гипотезы Хабибуллина. 
\medskip

\noindent{\bf Ключевые слова:} {гипотеза Хабибуллина, контрпример Шарипова, 
интегральные неравенства, интегральные преобразования.}
\medskip 

\end{quote}
}
\section{Введение.}\label{sec:intro}
    Гипотеза Хабибуллина --- это некоторое утверждение относительно
интегральных неравенств. Первоначально она была сформулирована в работах 
\cite{1} и \cite{2} в классе неотрицательных неубывающих функций, выпуклых 
относительно логарифма. В работе \cite{3} были даны две другие формулировки 
этой гипотезы, эквивалентные первоначальной. При этом условие выпуклости
относительно логарифма было снято и гипотеза Хабибуллина была сначала
переформулирована в классе неотрицательных неубывающих функций, а затем 
в классе неотрицательных непрерывных функций. Приведём все три формулировки 
гипотезы Хабибуллина. Первая из них --- это исходная формулировка. 
\begin{conjecture}[Хабибуллин]\label{conj:Khab1}
Пусть $\lambda$ --- положительное вещественное число и пусть $n$ --- 
целое число, такое, что $n\geqslant 2$. Пусть $S(x)$ --- 
неотрицательная неубывающая функция на интервале на $[0,+\infty)$, 
выпуклая относительно логарифма. Тогда, если неравенство 
\begin{equation*}
\hskip -2em
\int\limits^{\,1}_0\!S(tx)\,(1-x^2)^{n-2}
\,x\,dx\leqslant t^{\lambda}
\end{equation*}
выполняется для всех \ $0\leqslant t<+\infty$, то из него вытекает неравенство
\begin{equation*}
\hskip -2em
\int\limits^{+\infty}_0\!\!S(t)\,\frac{t^{2\lambda-1}}{(1+t^{2\lambda})^2}
\,dt\leqslant\frac{\pi\,(n-1)}{2\lambda}\,\prod^{n-1}_{k=1}\Bigl(1
+\frac{\lambda}{2k}\Bigr).
\end{equation*}
\end{conjecture}
     Условие выпуклости функции $S(x)$ относительно логарифма в гипотезе
\ref{conj:Khab1} означает, что функция $\sigma(x)=S(e^x)$ является выпуклой 
функцией в обычном понимании, т.\,е. при при любых значениях $x$ и $y$ она 
удовлетворяет неравенству
$$
\sigma(\alpha\,x+\beta\,y)\leqslant\alpha\,\sigma(x)+\beta\,\sigma(y)
\text{, \ если \ }\alpha\geqslant 0,\ \beta\geqslant 0\text{\ \ и \ }
\alpha+\beta=1.
$$
\begin{conjecture}[Хабибуллин]\label{conj:Khab2} Пусть $\alpha$ --- 
положительное вещественное число и пусть $n$ --- целое число, такое, 
что $n\geqslant 1$. Пусть $h(x)$ --- неотрицательная неубывающая функция 
на интервале $[0,+\infty)$. Тогда, если неравенство 
\begin{equation}\label{eq:1.1}
\hskip -2em
\int\limits^{\,1}_0 \frac{h(tx)}{x}\,(1-x)^{n-1}
\,dx\leqslant t^{\alpha}
\end{equation}
выполняется для всех \ $0\leqslant t<+\infty$, то из него вытекает 
неравенство
\begin{equation}\label{eq:1.2}
\hskip -2em
\int\limits^{+\infty}_0\frac{h(t)}{t} \,\frac{dt}{1+t^{2\alpha}}\leqslant
\frac{\pi}{2}\,\prod^{n-1}_{k=1}\Bigl(1+\frac{\alpha}{k}\Bigr).
\end{equation}
\end{conjecture}
\begin{conjecture}[Хабибуллин]\label{conj:Khab3} Пусть $\alpha$ --- 
положительное вещественное число и пусть $n$ --- натуральное число, такое
что $n\geqslant 1$. Пусть $q(x)$ --- неотрицательная непрерывная функция
на интервале $[0,+\infty)$. Тогда, если неравенство 
\begin{equation*}
\hskip -2em
\int\limits^{\,1}_0\left(\,\,\int\limits^{\,1}_x(1-y)^{n-1}
\,\frac{dy}{y}\right)q(tx)\,dx\leqslant t^{\alpha-1}
\end{equation*}
выполняется для всех \ $0\leqslant t<+\infty$, то из него вытекает 
неравенство
\begin{equation*}
\hskip -2em
\int\limits^{+\infty}_0\!\!q(t)\,\ln\Bigl(1+\frac{1}
{t^{\,2\kern 0.2pt\alpha}}\Bigr)\,dt\leqslant
\pi\,\alpha\prod^{n-1}_{k=1}\Bigl(1+\frac{\alpha}{k}\Bigr).
\end{equation*}
\end{conjecture}
     В работе \cite{4} было доказано, что гипотеза \ref{conj:Khab1} 
верна в случае $0<\lambda\leqslant 1$. Для гипотез \ref{conj:Khab2}
и \ref{conj:Khab3} это означает, что они верны при $0<\alpha\leqslant 1/2$.
\par
     В серии работ \cite{5,6,7} отдельно исследовалась третья формулировка
гипотезы Хабибуллина. При этом ещё раз было доказано, что гипотеза 
\ref{conj:Khab3} верна при $0<\alpha\leqslant 1/2$ для любых целых 
$n\geqslant 1$. Однако, как оказалось, за пределами интервала 
$0<\alpha\leqslant 1/2$ гипотеза \ref{conj:Khab3} не всегда верна. В
работе \cite{7} был построен контрпример к гипотезе \ref{conj:Khab3} 
для случая $n=2$ и $\alpha=2$. Основная цель данной работы пересчитать
построенный Шариповым в \cite{7} контрпример на случай гипотез
\ref{conj:Khab2} и \ref{conj:Khab1} соответственно.\par 
\section{Контрпример Шарипова.}\label{sec:shrcounterexample}
     Контрпример Шарипова к гипотезе \ref{conj:Khab3}, построенный
в работе \cite{7}, --- это функция $q(x)$, заданная на интервале 
$[0,+\infty)$ при помощи формулы
\begin{equation}\label{eq:2.1}
\hskip -2em
q(x)=\begin{cases} 12\,x\,(1-\varepsilon\,r(x))&\text{при \ }
0\leqslant x<x_0,\\
\kern 2.5em 12\,x&\text{при \ }x\geqslant x_0,
\end{cases}
\end{equation}
в которой $x_0=\root 4\of{3/5}$, а $r(x)$ --- некоторый конкретный полином
четвёртой степени по $x$. Контрпример \eqref{eq:2.1} содержит числовой
параметр $\varepsilon$, который может быть любым фиксированным вещественным 
числом, удовлетворяющим неравенству $0<\varepsilon\leqslant 1$. То есть
формула \eqref{eq:2.1} даёт не один контрпример, а однопараметрическое 
семейство контрпримеров, зависящее от параметра $\varepsilon$.\par
      Полином $r(x)$ в формуле \eqref{eq:2.1} выражается через другой полином 
четвёртой степени $R(\tau)$ при помощи линейной замены переменной:
\begin{equation}\label{eq:2.2}
\hskip -2em
r(x)=R(\tau)\text{, \ где \ }\tau=\frac{x_0-x}{x_0}. 
\end{equation}
А полином $R(\tau)$ в \eqref{eq:2.2} уже определяется явной формулой:
\begin{equation*}
\hskip -2em
R(\tau)=(21\,\tau^3-34\,\tau^2+16\,\tau-2)\,\tau. 
\end{equation*}\par
\section{Контрпример к гипотезе \ref{conj:Khab2}.}\label{sec:secondcounterexample}
     Контрпример к гипотезе \ref{conj:Khab2} должен быть неотрицательной
неубывающей функцией $h(x)$ на интервале $[0,+\infty)$, удовлетворяющей 
неравенству \eqref{eq:1.1}, но не удовлетворяющей неравенству \eqref{eq:1.2}.
При выводе гипотезы \ref{conj:Khab3} из гипотезы \ref{conj:Khab2} в работах
\cite{3} и \cite{4} была использована формула 
\begin{equation}\label{eq:3.1}
\hskip -2em
q(x)=\frac{dh(x)}{dx}.
\end{equation}
Согласно формуле \eqref{eq:3.1} для получения требуемой функции $h(x)$ надо
проинтегрировать функцию \eqref{eq:2.1}. Это даёт
\begin{equation}\label{eq:3.2}
\hskip -2em
h(x)=\!\!\int\limits^{\,x}_0\!q(y)\,dy+C,
\end{equation}
где $C$ --- константа интегрирования.\par
     Выбор константы интегрирования $C$ в формуле \eqref{eq:3.2} определяется 
неравенством \eqref{eq:1.1}. В этом неравенстве следует положить $\alpha=2$ и $n=2$, 
поскольку контрпример Шарипова \eqref{eq:2.1} относится именно к такому случаю. Тогда
неравенством \eqref{eq:1.1} запишется в следующем виде:
\begin{equation}\label{eq:3.3}
\hskip -2em
\int\limits^{\,1}_0 \frac{h(tx)}{x}\,(1-x)
\,dx\leqslant t^2.
\end{equation}
Подставим $t=0$ в неравенство \eqref{eq:3.3}. При этом аргумент функции $h(tx)$
занулится и её можно будет вынести за знак интегрирования:
\begin{equation}\label{eq:3.4}
\hskip -2em
h(0)\int\limits^{\,1}_0 \frac{1-x}{x}\,dx\leqslant 0.
\end{equation}
Интеграл в формуле \eqref{eq:3.4} положителен. По этой причине неравенство 
\eqref{eq:3.4} даёт $h(0)\leqslant 0$. С другой стороны, $h(x)$ должна быть
неотрицательной функцией на интервале $[0,+\infty)$. Отсюда $h(0)\geqslant 0$.
Из двух противоположных неравенств $h(0)\leqslant 0$ и $h(0)\geqslant 0$ 
вытекает следующее равенство:
\begin{equation}\label{eq:3.5}
\hskip -2em
h(0)=0.
\end{equation}
Равенство \eqref{eq:3.5} означает, что константа интегрирования $C$ в
формуле \eqref{eq:3.2} должна выбираться равной нулю, что даёт
\begin{equation}\label{eq:3.6}
\hskip -2em
h(x)=\!\!\int\limits^{\,x}_0\!q(y)\,dy.
\end{equation}\par
     Теперь функция $h(x)$ вычисляется прямой подстановкой функции 
\eqref{eq:2.1} в формулу \eqref{eq:3.6}. После вычисления интеграла
это даёт
\begin{equation}\label{eq:3.7}
\hskip -2em
h(x)=\begin{cases} 6\,x^2\,(1-\varepsilon\,u(x))&\text{при \ }
0\leqslant x<x_0,\\
\kern 2.5em 6\,x^2&\text{при \ }x\geqslant x_0.
\end{cases}
\end{equation}
Многочлен $u(x)$ в формуле \eqref{eq:3.7} определяется равенством
\begin{equation}\label{eq:3.8}
\hskip -2em
u(x)=U(\tau)\text{, \ где \ }\tau=\frac{x_0-x}{x_0}. 
\end{equation}
Через $U(\tau)$ в формуле \eqref{eq:3.8} обозначен многочлен четвёртой 
степени по переменной $\tau$, определяемый следующей явной формулой:
\begin{equation}\label{eq:3.9}
\hskip -2em
U(\tau)=(7\,\tau^2-8\,\tau+2)\,\tau^2. 
\end{equation}\par
\begin{theorem}\label{theo:3.1} Для каждого конкретного значения параметра
$\varepsilon$, удовлетворяющего неравенствам $0<\varepsilon\leqslant 1$,
функция $h(x)$, определяемая формулами \eqref{eq:3.7}, \eqref{eq:3.8} и
\eqref{eq:3.9}, является контрпримером к гипотезе Хабибуллина
\ref{conj:Khab2} для случая $n=2$ и $\alpha=2$. 
\end{theorem}
     Доказательство теоремы \ref{theo:3.1} может быть получено прямыми
вычислениями. Сама же функция \eqref{eq:3.7} с многочленом $u(x)$ из 
\eqref{eq:3.8} --- это есть контрпример Шарипова, пересчитанный со 
случая гипотезы \ref{conj:Khab3} на случай гипотезы \ref{conj:Khab2}.
\par
\section{Контрпример к гипотезе 
\ref{conj:Khab1}.}\label{sec:thirdcounterexample}
     Для того, чтобы пересчитать контрпример \eqref{eq:3.7} со случая
гипотезы \ref{conj:Khab2} на случай гипотезы \ref{conj:Khab1}, надо пройти
путем, проложенным в работе \cite{3}, но в обратном направлении. Сначала 
определим функцию
\begin{equation}\label{eq:4.1}
\hskip -2em
s(x)=4\,h(x^2).
\end{equation}
Подобно $h(x)$, функция $s(x)$ из \eqref{eq:4.1} является неотрицательной 
неубывающей функцией на интервале $[0,+\infty)$. С помощью $s(x)$ определим
функцию
\begin{equation}\label{eq:4.2}
\hskip -2em
S(x)=\!\!\int\limits^{\,\,x}_0\frac{s(t)}{t}\,dt
=\!\!\int\limits^{\,\,x}_0\frac{4\,h(t^2)}{t}\,dt.
\end{equation}
Преобразование \eqref{eq:4.2} основано на утверждении 5.1 из работы \cite{8},
которое было использовано в работе \cite{3}. Полученная в результате такого 
преобразования функция $S(x)$ будет неотрицательной, неубывающей и одновременно
выпуклой относительно логарифма функцией.\par
     Обозначим $x_1=\sqrt{x_0}=\root 8\of{3/5}$. После этого подставим функцию
\eqref{eq:3.7} в формулу \eqref{eq:4.2}. В результате прямых вычислений получим
\begin{equation}\label{eq:4.3}
\hskip -2em
S(x)=\begin{cases} 6\,x^4\,(1-\varepsilon\,v(x))&\text{при \ }
0\leqslant x<x_1,\\
\kern 2.5em 6\,x^4&\text{при \ }x\geqslant x_1.
\end{cases}
\end{equation}
Здесь $v(x)$ --- некоторый конкретный полином восьмой степени. Полином $v(x)$ 
в формуле \eqref{eq:4.3} выражается через другой полином восьмой степени 
$V(\theta)$ при помощи линейной замены переменной:
\begin{equation}\label{eq:4.4}
\hskip -2em
v(x)=V(\theta)\text{, \ где \ }\theta=\frac{x}{x_1}. 
\end{equation}
А полином $V(\theta)$ в \eqref{eq:4.4} уже определяется явной формулой:
\begin{equation}\label{eq:4.5}
\hskip -2em
V(\theta)=\frac{(7\,\theta^{\kern 0.3pt 2}-3)
\,(\theta^{\kern 0.3pt 2}-1)^3}{3}. 
\end{equation}\par 
Согласно результатам работы \cite{3} при переходе от гипотезы \ref{conj:Khab1}
к гипотезе \ref{conj:Khab2} параметр $n$ не меняется, а параметр $\lambda$
заменяется параметром $\alpha=\lambda/2$. По этой причине можно сформулировать
следующий результат. 
\begin{theorem}\label{theo:4.1} Для каждого конкретного значения параметра
$\varepsilon$, удовлетворяющего неравенствам $0<\varepsilon\leqslant 1$,
функция $S(x)$, определяемая формулами \eqref{eq:4.3}, \eqref{eq:4.4} и
\eqref{eq:4.5}, является контрпримером к гипотезе Хабибуллина \ref{conj:Khab1} 
для случая $n=2$ и $\lambda=4$. 
\end{theorem}Теорема \ref{theo:4.1} доказывается прямыми вычислениями. Функция $S(x)$, 
упоминаемая в этой теореме, --- это контрпример Шарипова, пересчитанный со 
уже со случая гипотезы \ref{conj:Khab2} на случай гипотезы \ref{conj:Khab1}.
\par
\vphantom{a}\vskip -25pt plus 25pt minus 0pt\vphantom{b}\par

\end{document}